\documentclass[12pt,amsmath]{article}
\usepackage{amsfonts}
\begin{document}

\newsymbol\rtimes 226F
\newsymbol\ltimes 226E
\newcommand{\text}[1]{\mbox{{\rm #1}}}
\newcommand{\Rep}{\text{Rep}}
\newcommand{\gr}{\text{gr}}
\newcommand{\Fun}{\text{Fun}}
\newcommand{\Hom}{\text{Hom}}
\newcommand{\End}{\text{End}}
\newcommand{\FPdim}{\text{FPdim}}
\newcommand{\GL}{\text{GL}}
\newcommand{\Sp}{\text{Sp}}
\newcommand{\Ps}{\text{Ps}}
\newcommand{\Ad}{\text{Ad}}
\newcommand{\ASp}{\text{ASp}}
\newcommand{\APs}{\text{APs}}
\newcommand{\Rad}{\text{Rad}}
\newcommand{\Corad}{\text{Corad}}
\newcommand{\SuperVect}{\text{SuperVect}}
\newcommand{\Vect}{\text{Vect}}
\newcommand{\Spec}{\text{Spec}}
\newcommand{\tr}{\text{tr}}
\newcommand{\cH}{{\cal H}}
\newcommand{\cR}{{\cal R}}
\newcommand{\cJ}{\cal J}
\newcommand{\gd}{\delta}
\newcommand{\lan}{\langle}
\newcommand{\ran}{\rangle}
\newcommand{\itms}[1]{\item[[#1]]}
\newcommand{\nin}{\in\!\!\!\!\!/}
\newcommand{\g}{{\bf g}}
\newcommand{\sub}{\subset}
\newcommand{\cntd}{\subseteq}
\newcommand{\go}{\omega}
\newcommand{\Pa}{P_{a^\nu,1}(U)}
\newcommand{\fx}{f(x)}
\newcommand{\fy}{f(y)}
\newcommand{\gD}{\Delta}
\newcommand{\gl}{\lambda}
\newcommand{\gL}{\Lambda}
\newcommand{\half}{\frac{1}{2}}
\newcommand{\sto}[1]{#1^{(1)}}
\newcommand{\stt}[1]{#1^{(2)}}
\newcommand{\Z}{\hbox{\sf Z\kern-0.720em\hbox{ Z}}}
\newcommand{\singcolb}[2]{\left(\begin{array}{c}#1\\#2
\end{array}\right)}
\newcommand{\ga}{\alpha}
\newcommand{\gb}{\beta}
\newcommand{\gga}{\gamma}
\newcommand{\ul}{\underline}
\newcommand{\ol}{\overline}
\newcommand{\qed}{\kern 5pt\vrule height8pt width6.5pt depth2pt}
\newcommand{\Lrraro}{\Longrightarrow}
\newcommand{\Nb}{|\!\!/}
\newcommand{\NN}{{\rm I\!N}}
\newcommand{\bsl}{\backslash}
\newcommand{\gt}{\theta}
\newcommand{\op}{\oplus}
\newcommand{\C}{{\bf C}}
\newcommand{\Q}{{\bf Q}}
\newcommand{\Op}{\bigoplus}
\newcommand{\CR}{{\cal R}}
\newcommand{\grr}{\omega_1}
\newcommand{\ben}{\begin{enumerate}}
\newcommand{\een}{\end{enumerate}}
\newcommand{\ndiv}{\not\mid}
\newcommand{\bab}{\bowtie}
\newcommand{\hal}{\leftharpoonup}
\newcommand{\har}{\rightharpoonup}
\newcommand{\ot}{\otimes}
\newcommand{\OT}{\bigotimes}
\newcommand{\bwe}{\bigwedge}
\newcommand{\eps}{\varepsilon}
\newcommand{\gs}{\sigma}
\newcommand{\rbraces}[1]{\left( #1 \right)}
\newcommand{\bbox}{$\;\;\rule{2mm}{2mm}$}
\newcommand{\sbraces}[1]{\left[ #1 \right]}
\newcommand{\bbraces}[1]{\left\{ #1 \right\}}
\newcommand{\OO}{_{(1)}}
\newcommand{\TT}{_{(2)}}
\newcommand{\FF}{_{(3)}}
\newcommand{\minus}{^{-1}}
\newcommand{\CV}{\cal V}
\newcommand{\CVs}{\cal{V}_s}
\newcommand{\un}{U_q(sl_n)'}
\newcommand{\on}{O_q(SL_n)'}
\newcommand{\slq}{U_q(sl_2)}
\newcommand{\olq}{O_q(SL_2)}
\newcommand{\UU}{U_{(N,\nu,\go)}}
\newcommand{\HH}{{\mathcal H}}
\newcommand{\ct}{\centerline}
\newcommand{\bs}{\bigskip}
\newcommand{\qua}{\rm quasitriangular}
\newcommand{\ms}{\medskip}
\newcommand{\noin}{\noindent}
\newcommand{\mat}[1]{$\;{#1}\;$}
\newcommand{\raro}{\rightarrow}
\newcommand{\map}[3]{{#1}\::\:{#2}\raro{#3}}
\newcommand{\alg}{{\rm Alg}}
\def\newtheorems{\newtheorem{theorem}{Theorem}[subsection]
                 \newtheorem{cor}[theorem]{Corollary}
                 \newtheorem{proposition}[theorem]{Proposition}
                 \newtheorem{lemma}[theorem]{Lemma}
                 \newtheorem{defn}[theorem]{Definition}
                 \newtheorem{Theorem}{Theorem}[section]
                 \newtheorem{Corollary}[Theorem]{Corollary}
                 \newtheorem{corollary}[theorem]{Corollary}
                 \newtheorem{Proposition}[Theorem]{Proposition}
                 \newtheorem{Lemma}[Theorem]{Lemma}
                 \newtheorem{Definition}[Theorem]{Definition}
                 \newtheorem{Example}[Theorem]{Example}
                 \newtheorem{Remark}[Theorem]{Remark}
                 \newtheorem{claim}[theorem]{Claim}
                 \newtheorem{sublemma}[theorem]{Sublemma}
                 \newtheorem{example}[theorem]{Example}
                 \newtheorem{definition}[theorem]{Definition}
                 \newtheorem{remark}[theorem]{Remark}
                 \newtheorem{question}[theorem]{Question}
                 \newtheorem{Question}[Theorem]{Question}
                 \newtheorem{Conjecture}[Theorem]{Conjecture}}
\newtheorems
\newcommand{\proof}{\par\noindent{\bf Proof:}\quad}
\newcommand{\dmatr}[2]{\left(\begin{array}{c}{#1}\\
                            {#2}\end{array}\right)}
\newcommand{\doubcolb}[4]{\left(\begin{array}{cc}#1&#2\\
#3&#4\end{array}\right)}
\newcommand{\qmatrl}[4]{\left(\begin{array}{ll}{#1}&{#2}\\
                            {#3}&{#4}\end{array}\right)}
\newcommand{\qmatrc}[4]{\left(\begin{array}{cc}{#1}&{#2}\\
                            {#3}&{#4}\end{array}\right)}
\newcommand{\qmatrr}[4]{\left(\begin{array}{rr}{#1}&{#2}\\
                            {#3}&{#4}\end{array}\right)}
\newcommand{\smatr}[2]{\left(\begin{array}{c}{#1}\\
                            \vdots\\{#2}\end{array}\right)}

\newcommand{\ddet}[2]{\left[\begin{array}{c}{#1}\\
                           {#2}\end{array}\right]}
\newcommand{\qdetl}[4]{\left[\begin{array}{ll}{#1}&{#2}\\
                           {#3}&{#4}\end{array}\right]}
\newcommand{\qdetc}[4]{\left[\begin{array}{cc}{#1}&{#2}\\
                           {#3}&{#4}\end{array}\right]}
\newcommand{\qdetr}[4]{\left[\begin{array}{rr}{#1}&{#2}\\
                           {#3}&{#4}\end{array}\right]}

\newcommand{\qbracl}[4]{\left\{\begin{array}{ll}{#1}&{#2}\\
                           {#3}&{#4}\end{array}\right.}
\newcommand{\qbracr}[4]{\left.\begin{array}{ll}{#1}&{#2}\\
                           {#3}&{#4}\end{array}\right\}}

\title{{\bf Triangular Hopf Algebras
With The Chevalley Property}}
\author{
Nicol\'as Andruskiewitsch $^1$
\and
Pavel Etingof $^2$
\and
Shlomo Gelaki $^3$
}
\footnotetext[1]
{FaMAF - UNC (5000) Ciudad,
Universitaria C\'ordoba,
Argentina,
{\rm email: andrus@mate.uncor.edu}
}
\footnotetext[2]
{MIT,
Department of Mathematics,
77 Massachusetts Avenue,
Cambridge, MA 02139, USA, \&
Columbia University,
Department of Mathematics,
2990 Broadway,
New York, NY 10027, USA,
{\rm email: etingof@math.mit.edu}
}
\footnotetext[3]
{Technion-Israel Institute of Technology,
Department of Mathematics,
Haifa 32000, Israel,
{\rm email: gelaki@math.technion.ac.il}
}
\maketitle

\section{Introduction}
Triangular Hopf algebras were introduced by Drinfeld [Dr].
They are the Hopf algebras whose representations
form a symmetric tensor category. In that sense, they are
the class of Hopf algebras closest to group algebras.
The structure of triangular Hopf algebras is far from trivial,
and yet is more tractable than that of general Hopf
algebras, due to their proximity to groups and Lie algebras.
This makes triangular Hopf algebras an excellent testing
ground for
general Hopf algebraic ideas, methods and conjectures.

A general classification
of triangular Hopf algebras is not known yet.
However, there are two classes
that are relatively well understood.
One of them is semisimple triangular Hopf algebras over $\C$,
for which a complete classification is given in [EG1,EG2].
The key theorem about such Hopf algebras states that each
of them is
obtained by twisting a group algebra of a finite group
(see [EG1, Theorem 2.1]). The proof of
this theorem is based on Deligne's theorem on Tannakian
categories [De2].

Another important class of Hopf algebras is that of
{\em pointed} ones. These are Hopf algebras whose all
simple comodules are $1-$dimensional.
Theorem 5.1 in [G] gives a classification of minimal
triangular pointed
Hopf algebras (we note that
the additional assumption made in [G, Theorem 5.1] is
superfluous, by Theorem \ref{gen} below).

Recall that a finite-dimensional algebra is called {\it
basic} if all of its simple modules are
$1-$dimensional
(i.e. if its dual is a pointed coalgebra).
The same Theorem 5.1 of [G] gives a classification
of minimal triangular basic Hopf algebras, since the dual of a
minimal
triangular Hopf algebra is again minimal triangular.

Basic and semisimple
Hopf algebras share a common property. Namely, the
Jacobson radical $\Rad(H)$ of such a Hopf algebra $H$ is a
Hopf ideal,
and hence the quotient $H/\Rad(H)$ (the semisimple part)
is itself a Hopf algebra. The representation-theoretic
formulation
of this property is: The tensor product of two simple
$H$-modules
is semisimple. A remarkable classical theorem of Chevalley
[C, p.88] states that, over $\C,$ this property holds for
the group algebra
of any (not necessarily finite) group. So let us call
this
property of $H$ {\bf the Chevalley property}.

The Chevalley property certainly fails for many
finite-dimensional Hopf algebras, e.g. for Lusztig's
finite-dimensional
quantum groups $\text{U}_q({\bf g})'$ at roots of unity
(also known as Frobenius-Lusztig kernels) [L].
However, we found that this property holds
for all examples we know of finite-dimensional {\it triangular}
Hopf algebras in characteristic $0.$ We felt, therefore, that it is
natural to
classify all finite-dimensional triangular Hopf algebras
with the Chevalley property. This is what we do in this paper.

We start by classifying triangular Hopf
algebras with $R-$matrix of rank $\le 2$.
We show that such a Hopf
algebra
is a suitable modification of a cocommutative Hopf
superalgebra
(i.e. the group algebra of a supergroup). On the other
hand, by a theorem of Kostant [Ko],
a finite supergroup is a semidirect
product of
a finite group with an odd vector space on which this group
acts.

Next we prove our main result:
Any finite-dimensional triangular Hopf algebra
with the Chevalley property
is obtained by twisting a triangular
Hopf algebra with $R-$matrix of rank $\le 2$. We also prove the
converse result
that any such Hopf algebra does have the Chevalley property.
As a corollary, we prove that any finite-dimensional triangular
Hopf algebra whose coradical is a Hopf subalgebra
(e.g. pointed) is obtained by twisting a triangular
Hopf algebra with $R-$matrix of rank $\le 2$.

The structure of the paper is as follows.

In Section 2
we give the definitions of Hopf superalgebras and twists
for them. We also discuss cocommutative Hopf
superalgebras, and
describe their classification (Kostant's theorem [Ko]).

In Section 3 we establish a correspondence
between usual Hopf algebras
and Hopf superalgebras, and show how this correspondence
extends
to twists and to triangular Hopf algebras.

In Section 4 we discuss the Chevalley property.

In Section 5 we prove our main result,
and discuss its consequences and some open questions.

In Section 6, using the main theorem,
we show that a finite-dimensional cotriangular pointed
Hopf
algebra is generated by its grouplike and skewprimitive
elements.
Thus we confirm the conjecture that this is the case
for any finite-dimensional pointed Hopf algebra over $\C$
[AS2],
in the cotriangular case. This allows us to strengthen the
main result of [G].

In Section 7 we prove that the categorical
dimensions of objects in any abelian symmetric rigid category
with finitely many irreducible objects are integers.
In particular, this is the case for the representation category
of a triangular Hopf algebra. This gives supporting evidence
for a positive answer to the question we ask in Section 5:
Is any finite-dimensional triangular Hopf algebra a twist of a
modified supergroup algebra?

In the appendix we use the Lifting method [AS1,AS2] to
give other proofs of Theorem \ref {1} and Corollary \ref
{mtp}, and a generalization of Lemma \ref {5}.

We note that similarly to the case of semisimple Hopf
algebras, the proof of our main result is based on
Deligne's theorem [De2]. In fact, we use Theorem 2.1 of
[EG1] to prove the main result of this paper.

\noin
{\bf Acknowledgments.}
The first author was partially supported by CONICET, CONICOR and
Secyt (UNC), Argentina.
The research of the second author was conducted for the
Clay Mathematical Institute. The first and third authors
are grateful to MSRI for its support. The third author is
grateful to Bar Ilan University for its support.

\section{Hopf Superalgebras}

\subsection{Supervector Spaces}

The ground field in this paper will always be the
field $\C$ of complex numbers.

We start by recalling the definition of the category of
supervector spaces.
A Hopf algebraic way to define this category is as follows.

Let $u$ be the generator of the group $\Z_2$ of two
elements, and set
\begin{equation}\label{u}
R_u:=\frac{1}{2}
(1\ot 1+1\ot u+u\ot 1-u\ot u)\in \C[\Z_2]\otimes \C[\Z_2].
\end{equation}
Then $(\C[\Z_2],R_u)$ is a
minimal triangular Hopf algebra.

\begin{definition} The category of supervector spaces
over $\C$
is the symmetric tensor category $\Rep(\C[\Z_2],R_u)$
of representations of the triangular Hopf algebra
$(\C[\Z_2],R_u).$ This category will be denoted by
$\SuperVect.$
\end{definition}

For $V\in \SuperVect$ and $v\in V$, we say that $v$ is
even
if $uv=v$ and odd if $uv=-v$. The set of even vectors in
$V$ is
denoted by $V_0$ and the set of odd vectors by $V_1$, so
$V=V_0\oplus V_1$.
We define the parity of a vector $v$ to be $p(v)=0$ if $v$ is
even and $p(v)=1$ if $v$ is odd (if $v$ is neither odd
nor even,
$p(v)$ is not defined).

Thus, as an ordinary tensor category, SuperVect is equivalent
to the category of representations of $\Z_2$, but the
commutativity
constraint is different
from that of $\Rep(\Z_2)$ and equals $\beta:=R_uP$, where
$P$ is the
permutation of components.
In other words, we have
\begin{equation}\label{symm}
\beta (v\otimes w)=(-1)^{p(v)p(w)}w\otimes v,
\end{equation}
where both $v,w$ are either even or odd.

\subsection{Hopf Superalgebras}

Recall that in any symmetric (more generally,
braided) tensor category,
one can define an algebra, coalgebra,
bialgebra, Hopf algebra, triangular Hopf algebra, etc,
to be an object of this category equipped with
the usual structure maps
(morphisms in this category),
subject to the same axioms as in the usual case.
In particular, any of these algebraic structures
in the category SuperVect is usually identified by the prefix
``super''. For example:

\begin{definition}
A Hopf superalgebra is a Hopf algebra
in
$\SuperVect.$
\end{definition}

More specifically, a Hopf superalgebra $\cH$
is an ordinary
$\Z_2$-graded associative unital algebra with
multiplication $m$,
equipped with a coassociative map
$\Delta:\cH\raro
\cH\ot \cH$ (a morphism in $\SuperVect$) which is
multiplicative in
the super-sense, and with a
counit and antipode satisfying the standard axioms.
Here multiplicativity in the super-sense
means that $\Delta$ satisfies the relation
\begin{equation}\label{sd}
\Delta(ab)=\sum (-1)^{p(a_2)p(b_1)}a_1b_1\ot a_2b_2
\end{equation}
for all $a,b\in \cH$
(where $\Delta(a)=\sum a_1\ot a_2$, $\Delta(b)=\sum
b_1\ot b_2$).
This is because the tensor product of
two algebras $A,B$ in $\SuperVect$ is defined to be $A\ot
B$ as a vector space, with multiplication
\begin{equation}\label{sm}
(a\ot b)(a'\ot b'):=(-1)^{p(a')p(b)}a a'\ot
bb'.
\end{equation}

\begin{remark} {\rm Hopf superalgebras appear in
[Ko], under the name of ``graded Hopf algebras''.
}
\end{remark}

Similarly, a
(quasi)triangular Hopf superalgebra $(\cH,\cR)$ is a Hopf
superalgebra
with an $R-$matrix (an even element $\cR\in \cH\otimes
\cH$)
satisfying the usual axioms. As in the even case, an
important role is played by the Drinfeld element $u$ of
$(\cH,\cR)$:
\begin{equation}\label{du}
u:=m\circ \beta \circ (Id\ot S)(\cR).
\end{equation}
For instance, $(\cH,\cR)$ is triangular if and only if $u$
is a grouplike element of $\cH.$

As in the even case, the tensorands of the $R-$matrix
of a (quasi)triangular Hopf superalgebra $\cH$ generate a
finite-dimensional
sub Hopf superalgebra ${\cH}_m$, called the {\em minimal
part of} $\cH$ (the proof does not differ essentially
from the proof of the analogous fact for Hopf algebras).
A (quasi)triangular Hopf superalgebra is said to be minimal if
it coincides with its minimal part. The dimension of the minimal part
is the {\it rank} of the $R-$matrix.

\subsection{Cocommutative Hopf Superalgebras}

\begin{definition}
We will say that a Hopf superalgebra $\cH$ is
commutative (resp. cocommutative) if $m=m\circ
\beta$
(resp. $\Delta=\beta\circ \Delta$).
\end{definition}

\begin{example}\label{supergr}
{\rm {\bf [Ko]} Let $G$ be a group, and $\g$ a Lie
superalgebra with an action of
$G$ by automorphisms of Lie superalgebras.
Let ${\cal H}:=\C[G] \ltimes \text{U}(\g),$ where
$\text{U}(\g)$
denotes the universal
enveloping algebra of $\g$.
Then $\cH$ is a
cocommutative Hopf
superalgebra, with
$\Delta(x)=x\ot 1+1\ot x,$ $x\in \g$, and $\Delta(g)=g\ot g,$
$g\in G$. In this Hopf superalgebra, we have
$S(g)=g^{-1}$, $S(x)=-x$, and in particular $S^2=Id$.

The Hopf superalgebra $\cH$ is finite-dimensional if and
only if
$G$ is finite, and $\g$ is finite-dimensional and purely
odd (and
hence commutative).
Then
${\cal H}=\C[G]\ltimes \Lambda V$, where $V=\g$ is an odd
vector space
with a $G$-action. In this case, $\cH^*$ is a commutative
Hopf
superalgebra.
}
\end{example}

\begin{remark} {\rm We note that as in the even case, it
is convenient to think about $\cH$ and $\cH^*$ in geometric terms.
Consider, for instance, the finite-dimensional case. In this case,
it is useful to think of the ``affine algebraic supergroup''
$\tilde G:=G\ltimes V$. Then one can regard $\cH$ as the group
algebra $\C[\tilde G]$ of this supergroup, and $\cH^*$ as its
function algebra $F(\tilde G)$. Having this in mind, we will call
the algebra $\cH$ {\bf a supergroup algebra}. }
\end{remark}

It turns out that like in the even case,
any cocommutative Hopf superalgebra is of the type
described
in Example \ref{supergr}. Namely, we have the following
theorem.

\begin{theorem}\label{kostant} {\bf ([Ko], Theorem 3.3)}
Let $\cH$ be a cocommutative Hopf superalgebra over $\C.$
Then $\cH=\C[{\bf G}(\cH)]\ltimes
\text{U}({\bf P}(\cH)),$
where $\text{U}({\bf P}(\cH))$ is the universal enveloping
algebra
of the
Lie superalgebra of primitive elements of $\cH,$
and ${\bf G}(\cH)$ is the group of grouplike elements of
$\cH.$
\end{theorem}

In particular, in the finite-dimensional case we get:

\begin{corollary}\label{kostantf}
Let $\cH$ be a finite-dimensional cocommutative
Hopf superalgebra over $\C.$
Then ${\cal H}=\C[{\bf G}({\cal H})]\ltimes \Lambda V,$
where $V$ is the space of primitive elements of $\cH$
(regarded as an odd vector space) and ${\bf G}({\cal H})$ is the
finite group of grouplikes of ${\cal H}.$
In other words, $\cH$ is a supergroup algebra.
\qed
\end{corollary}

This corollary will be used below, and
although it follows at once from Theorem \ref{kostant},
for the sake of completeness we will give its proof in Section
5.

\subsection{Twists}

A twist for a Hopf algebra in any symmetric tensor category
is defined in the same way as in the usual case (see
[Dr]). However, for the
reader's convenience, we will repeat this definition
(for Hopf superalgebras).

Let $\cH$ be a Hopf superalgebra. The multiplication,
unit,
comultiplication, counit and antipode in $\cH$ will be
denoted by $m,1,\gD,\varepsilon,S$ respectively.

\begin{definition}\label{t}
A twist for $\cH$ is an invertible even
element ${\cal J}\in \cH\ot \cH$ which satisfies
\begin{equation}\label{t1}
(\Delta\ot Id)({\cal J})({\cal J}\ot 1)=(Id\ot \Delta)({\cal
J})(1\ot {\cal J})\;\;and \;\;
(\varepsilon\ot
Id)({\cal J})=(Id\ot \varepsilon)({\cal J})=1,
\end{equation}
where $Id$ is the identity map of $\cH.$
\end{definition}

Given a twist ${\cal J}$ for $\cH,$ one can define a new
Hopf superalgebra structure
$$({\cal H}^{{\cal J}},m,1,\Delta^{{\cal
J}},\varepsilon,S^{{\cal J}})$$
on the algebra $({\cal H},m,1)$
as follows. The coproduct is determined by
\begin{equation}\label{t2}
\Delta^{{\cal J}}(a)={{\cal J}}^{-1}\Delta(a){\cal J}\;\text{for
any}\;a\in
\cH,
\end{equation}
and the antipode is determined by
\begin{equation}\label{tant}
S^{{\cal J}}(a)=Q^{-1}S(a)Q\;\text{for any}\;a\in \cH,
\end{equation}
where $Q:=m\circ(S\ot Id)({\cal J}).$

If $\cH$ is (quasi)triangular with the universal
$R-$matrix $\cR$ then so
is ${\cal H}^{{\cal J}}$ with the universal $R-$matrix
${\cal R}^{{\cal J}}:={\cal J}_{21}^{-1}{\cal R} {\cal J}.$

\section{Triangular Hopf Algebras With Drinfeld
Element of Order $\le 2$}

\subsection{The Correspondence Between Hopf Algebras and
Superalgebras}

We can now prove our first results,
which will be
essential in the next section.
We start with a correspondence theorem
between Hopf algebras and Hopf superalgebras.

\begin{theorem}\label{cores0}
There is a one to one correspondence between:
\ben
\item isomorphism classes of pairs $(H,u)$ where $H$ is an
ordinary Hopf
algebra, and $u$ is a grouplike element in $H$ such that
$u^2=1,$ and
\item
isomorphism classes of pairs $(\cH,g)$ where $\cH$
is a Hopf superalgebra, and
$g$ is a grouplike element in $\cH$
such that $g^{2}=1$ and $g x
g^{-1}=(-1)^{p(x)}x$ (i.e. $g$ acts on $x$ by
its parity),
\een
such that the tensor categories of representations of $H$ and
$\cH$ are equivalent.
\end{theorem}

\proof Let $(H,u)$ be an ordinary
Hopf algebra
with comultiplication $\Delta,$ counit $\eps$,
antipode $S$,
and a grouplike element $u$
such that $u^2=1.$
Let $\cH=H$ regarded as a superalgebra, where the
$\Z_2-$grading is
given by the adjoint
action of $u.$ For $h\in H,$ let us define $\Delta_0,
\Delta_1$ by writing $\Delta(h)=\Delta_0(h)+\Delta_1(h),$
where
$\Delta_0(h)\in H\ot H_0$ and $\Delta_1(h)\in H\ot H_1.$
Define a map
$\tilde\Delta:\cH\raro \cH\ot \cH$ by
$\tilde\Delta(h):=\Delta_0(h)-(-1)^{p(h)}(u\ot
1)\Delta_1(h).$ Define
$\tilde S(h):=u^{p(h)}S(h),$ $h\in H.$
Then it is straightforward
to verify that $(\cH,\tilde\Delta,\varepsilon,\tilde S)$
is a Hopf superalgebra.

The element $u$ remains grouplike in the new Hopf superalgebra,
and acts by parity, so we can set $g:=u$.

Conversely, suppose that $(\cH,g)$ is a pair where
$\cH$ is a Hopf superalgebra with
comultiplication $\tilde \Delta,$ counit $\varepsilon$,
antipode $\tilde S$, and a grouplike element $g$, with $g^2=1$,
acting by parity.
For $h\in \cH,$ let us define
$\tilde {\Delta}_0,\tilde {\Delta}_1$ by writing
$\tilde{\Delta}(h)=\tilde{\Delta}_0(h)+\tilde{\Delta}_1(h),$
where $\tilde{\Delta}_0(h)\in \cH\ot \cH_0$ and
$\tilde{\Delta}_1(h)\in \cH\ot \cH_1.$
Let $H=\cH$ as algebras, and
define a map
$\Delta:H\raro H\ot H$ by
$\Delta(h):=\tilde{\Delta}_0(h)-(-1)^{p(h)}(
g\ot 1)\tilde{\Delta}_1(h).$
Define $S(h):=g^{p(h)}\tilde S(h),$ $h\in H.$ Then it is
straightforward to verify that $(H,\Delta,\varepsilon,S)$ is an
ordinary Hopf
algebra, and we can set $u:=g$.

It is obvious that the two assignments constructed
above are inverse to each other.
The equivalence of tensor categories
is straightforward to verify. The theorem is proved.\qed

Theorem \ref{cores0} implies the following. Let $\cH$ be {\em
any} Hopf superalgebra, and $\C[\Z_2]\ltimes \cH$ be the
semidirect product, where the generator $g$ of $\Z_2$ acts on
$\cH$ by $gxg^{-1}=(-1)^{p(x)}x$. Then we can define an ordinary
Hopf algebra $\ol {\cH}$, which is the one corresponding to
$(\C[\Z_2]\ltimes \cH,g)$ under the correspondence of Theorem
\ref{cores0}.

The constructions of this section have the following explanation
in terms of Radford's biproduct construction [R2]. Namely $\cH$ is
a Hopf algebra in the Yetter-Drinfeld category of $\C[\Z_2],$ so
Radford's biproduct construction yields a Hopf algebra structure
on $\C[\Z_2]\ot \cH,$ and it is straightforward to see that this
Hopf algebra is exactly $\ol {\cH}.$ Moreover, it is clear that
for any pair $(H,u)$ as in Theorem \ref{cores0}, $gu$ is central
in $\ol {\cH}$ and $H=\ol {\cH}/(gu-1).$

Let us give an interesting corollary
of Theorem \ref{cores0}, even though we will not use it.

\begin{corollary}\label{ss} Let $\cH$ be a
finite-dimensional Hopf
superalgebra over $\C$. Then:
\ben

\item $\cH$ is semisimple if and only if it is
cosemisimple.

\item If $\cH$ is semisimple then $S^4=Id$.

\item If $\cH$ is semisimple and $S^2=Id$, then $\cH$
is purely even, i.e. it is a usual semisimple Hopf algebra.
\een
\end{corollary}

\proof
1. If $\cH$ is semisimple then so is $\ol {\cH}$,
hence so is $(\ol {\cH})^*$.
But it is easy to show that $(\ol {\cH})^*$ is isomorphic
as
an algebra
to $\C[\Z_2]\ltimes \cH^*$ (unlike the dual of
$\C[\Z_2]\ltimes \cH$,
which is isomorphic to $\C[\Z_2]\otimes \cH^*$).
Thus this crossed product algebra is semisimple.
It is well known (and easy to show)
that this implies the semisimplicity of $\cH^*$.

2. The Hopf algebra $\ol {\cH}$ is semisimple,
so we have $S^2=Id$ in it. Thus, in $\cH$ we have
$S^2=\Ad(g)$, so $S^4=\Ad(g^2)=Id$.

3. Since $S^2=\Ad(g),$ $g$ has to be central. Thus, $\cH$
is purely even. \qed

\begin{remark} {\rm The example of supergroup algebras
shows that
for finite-dimensional Hopf superalgebras, unlike
usual Hopf algebras, $S^2=Id$ does not imply
semisimplicity or
cosemisimplicity. In fact, Corollary \ref{ss}, part 3,
shows that,
in a sense, the situation is exactly the opposite.
}
\end{remark}

\subsection{Correspondence of Twists}

Let us say that a twist $J$ for a Hopf algebra $H$
with an involutive grouplike element
$g$ is {\it even} if
it is invariant under $\Ad(g)$.

\begin{proposition}\label{1a}
Let $(\cH,g)$ be a pair as in Theorem
\ref{cores0}, and let
$H$ be the associated ordinary Hopf algebra. Let
${\cal J}\in \cH\ot \cH$ be an even element. Write
${\cal J}=
{\cal J}_0+{\cal J}_1,$ where ${\cal J}_0\in \cH_0\ot
\cH_0$ and
${\cal J}_1\in \cH_1\ot \cH_1.$ Define
$J:={\cal J}_0-
(g\ot 1) {\cal J}_1$. Then $J$ is an even twist for $H$
if and only if ${\cal J}$ is a twist for $\cH.$
Moreover, $\cH^{{\cal J}}$ corresponds to $H^J$ under the
correspondence in Theorem \ref{cores0}.
Thus, there is a one to one correspondence between even twists
for $H$ and twists for $\cH$, given by $J\to {\cal J}$.
\end{proposition}

\proof Straightforward.\qed

\subsection{The Correspondence Between Triangular Hopf
Algebras and Superalgebras}

Let us now return to our main subject, which is
triangular Hopf algebras and superalgebras.
For triangular Hopf algebras whose Drinfeld element $u$ is
involutive, we will make the natural choice of the
element $u$ in Theorem \ref{cores0}, namely define it to be the
Drinfeld element of $H$.

\begin{theorem}\label{cores}
The correspondence of Theorem \ref{cores0}
extends to a one to one correspondence between:
\ben
\item isomorphism classes of ordinary
triangular Hopf
algebras $H$ with Drinfeld element $u$ such that $u^2=1,$ and
\item
isomorphism classes of pairs $(\cH,g)$ where $\cH$
is a
triangular Hopf superalgebra with Drinfeld element $1$ and
$g$ is an element of ${\bf G}(\cH)$
such that $g^{2}=1$ and $gxg^{-1}=(-1)^{p(x)}x.$
\een
\end{theorem}

\proof Let $(H,R)$ be a triangular Hopf algebra with $u^2=1$.
Since $(S\ot S)(R)=R$ and $S^2=\Ad(u)$ [Dr], $u\ot u$ and $R$
commute. Hence we can write $R=R_0+R_1$, where $R_0\in H_0\ot
H_0$ and $R_1\in H_1\ot H_1.$ Let $\cR:=(R_0+(1\ot u)R_1)R_u$.
Then $\cR$ is even. Indeed, since $R_0=1/2(R+(u\ot 1)R(u\ot 1))$
and $R_1=1/2(R-(u\ot 1)R(u\ot 1)),$ $u\ot u$ and $\cR$ commute.

It is now straightforward to show that $(\cH,\cR)$ is triangular
with Drinfeld element $1.$ Let us show for instance that $\cR$ is
unitary. Let us use the notation $a*b, X^{21}$ for multiplication
and opposition in the tensor square of a superalgebra, and the
notation $ab,X^{op}$ for usual algebras. Then,
$$\cR*\cR_{21}=
(R_0+(1\ot u)R_1)R_u*(R_0^{op} - (u\ot 1)R_1^{op})R_u.$$ Since,
$R_uR_0=R_0R_u,$ $R_uR_1=-(u\ot u)R_1R_u,$ we get that the RHS
equals
$$(R_0+(1\ot u)R_1)*(R_0^{op}+(1\ot u)R_1^{op})=
R_0R_0^{op}+R_1R_1^{op}+(1\ot u)(R_1R_0^{op}+R_0R_1^{op}).$$
But, $R_0R_0^{op}+R_1R_1^{op}=1$ and $(1\ot
u)(R_1R_0^{op}+R_0R_1^{op})=0,$ since $RR^{op}=1,$ so we
are done.

Conversely, suppose that $(\cH,g)$ is a pair where
$\cH$ is a triangular Hopf superalgebra with $R-$matrix
$\cR$
and Drinfeld element $1.$ Let
$\cR=\cR_0+\cR_1$, where
$\cR_0$ has even components, and $\cR_1$ has odd
components. Let $R:=(\cR_0+(1\ot g)\cR_1) R_g$.
Then it is straightforward
to show that $(H,R)$ is
triangular with Drinfeld element $u=g$.
The theorem is proved.\qed

\begin{corollary}\label{mintr}
If $(\cH,\cR)$ is a triangular
Hopf superalgebra with Drinfeld element 1,
then the Hopf algebra $\overline{\cH}$ is also
triangular, with the $R-$matrix
\begin{equation}\label{olr}
\ol R:=(\cR_0+(1\otimes g)\cR_1)R_g,
\end{equation}
where $g$ is the grouplike element adjoined to $\cH$ to
obtain $\overline{\cH}$. Moreover, $\cH$ is minimal if and
only if so is $\ol {\cH}.$
\end{corollary}

\proof Clear.\qed

The following corollary, combined with Kostant's theorem, gives a
classification of triangular Hopf algebras
with $R$-matrix of rank $\le 2$ (i.e. of the form
$R_u$ as in (\ref{u}), where $u$ is a grouplike of order
$\le 2$).

\begin{corollary}\label{cores'}
The correspondence of Theorem \ref{cores}
restricts to a one to one correspondence between:
\ben
\item isomorphism classes of ordinary
triangular Hopf algebras with $R-$matrix of
rank $\le 2,$ and
\item
isomorphism classes of pairs $(\cH,g)$ where
$\cH$ is a cocommutative Hopf superalgebra and $g$ is an
element of ${\bf G}(\cH)$
such that $g^{2}=1$ and $g x
g^{-1}=(-1)^{p(x)}x$.
\een
\end{corollary}
\proof Let $(H,R)$ be an ordinary triangular Hopf algebra with
$R-$ matrix of rank $\le 2.$ In particular, the Drinfeld element
$u$ of $H$ satisfies $u^2=1,$ and $R=R_u$. Hence by Theorem
\ref{cores}, $(\cH,\tilde\Delta,\cR)$ is a triangular Hopf
superalgebra. Moreover, it is cocommutative since $\cR=R_uR_u=1.$

Conversely, for any $(\cH,g)$,
by Theorem \ref{cores}, the pair $(H,R_{g})$ is an
ordinary triangular Hopf algebra, and clearly the
rank of $R_{g}$ is $\le 2.$ \qed

In particular, Corollaries \ref{kostantf} and \ref{cores'}
imply that finite-dimensional
triangular Hopf algebras with $R-$matrix of rank $\le 2$
correspond to supergroup algebras.
In view of this, we make the following definition.

\begin{definition}
A finite-dimensional triangular Hopf algebra with $R-$matrix of
rank $\le 2$ is called a modified supergroup algebra.
\end{definition}

\subsection{Construction of Twists for Supergroup
Algebras}

\begin{proposition}\label{st}
Let $\cH=\C[G] \ltimes \Lambda V$ be a supergroup algebra.
Let $r\in S^2V.$ Then ${\cal J}:=e^{r/2}$ is a twist for $\cH.$
Moreover,
$((\Lambda V)^{{\cal J}},{\cal J}_{21}^{-1}{\cal J})$ is minimal
triangular if
and only if $r$ is nondegenerate.
\end{proposition}
\proof Straightforward. \qed

\begin{example}
{\rm Let $G$ be the group of order $2$ with generator $g.$
Let $V:=\C$ be the nontrivial $1-$dimensional
representation of $G,$ and write
$\Lambda V=sp\{1,x\}.$ Then the associated ordinary
triangular Hopf algebra
to $(\cH,g):=(\C[G] \ltimes \Lambda V,g)$ is
Sweedler's $4-$dimensional
Hopf algebra $H$ [S] with the triangular structure $R_g.$
Namely, the algebra $H$ is generated by a grouplike
element $g$ and a $1:g$ skew
primitive element $x$ (i.e. $\Delta(x)=x\ot 1+g\ot x$)
satisfying
the relations $g^2=1,$ $x^2=0$ and $gx=-xg.$
It is known [R1] that
the set of triangular structures on $H$ is parameterized
by $\C;$ namely,
$R$ is a triangular structure on $H$ if and only if
$$R=R_{\lambda}:=
R_g-\frac{\lambda}{2}(x\ot x-gx\ot x+x\ot gx +gx\ot gx),
\;\lambda\in \C.$$ Clearly,
$(H,R_{\lambda})$ is minimal if and only if $\lambda\ne 0.$

Let $r\in S^2V$ be defined by $r:=\lambda x\otimes x,$
$\lambda\in \C.$ Set
${\cal J}_\lambda:=e^{r/2}=1+\frac{1}{2}\lambda x\ot x$;
it is a
twist for $\cH.$
Hence, $J_{\lambda}:=1-\frac{1}{2}\lambda gx\ot x$ is a
twist for $H.$ It is
easy to check that
$R_{\lambda}=(J_{\lambda})_{21}^{-1}R_gJ_{\lambda}.$
Thus, $(H,R_\lambda)=(H,R_0)^{J_\lambda}$.
}
\end{example}

\begin{remark}{\rm  In fact, Radford's classification
of triangular structures on $H$ can be easily deduced
from Lemma \ref{5} below.
}
\end{remark}

\section{The Chevalley Property}

Recall that in the introduction we made the following
definition.

\begin{Definition} A Hopf algebra $H$ over $\C$ is said to have
the Chevalley property if the tensor product of any two simple
$H$-modules is semisimple. More generally, let us say that a
tensor category has the Chevalley property if the tensor product
of two simple objects is semisimple.
\end{Definition}

Let us give some equivalent formulations of the Chevalley
property.

\begin{Proposition}\label{chev}
Let $H$ be a finite-dimensional Hopf algebra over $\C$
and let $A
:= H^*$. The following
conditions  are equivalent:
\ben
\item $H$ has the Chevalley property.
\item The category of (right) $A$-comodules has the
Chevalley property.
\item $\Corad(A)$ is a Hopf subalgebra of $A.$
\item $\Rad(H)$ is a Hopf ideal and thus $H/\Rad(H)$
is a Hopf algebra.
\item $S^2=Id$ on $H/\Rad(H)$, or equivalently on $\Corad(A)$.
\een
\end{Proposition}

\proof
(1. $\Leftrightarrow$ 2.) Clear, since the categories of left
$H$-modules
and right $A$-comodules are equivalent.

(2. $\Rightarrow$ 3.)
Recall the definition of a matrix coefficient of a comodule $V$
over $A$. If $\rho: V \to V \otimes A$ is the coaction, $v\in
V$, $\alpha
\in V^*$,
then
$$
\phi^V_{v, \alpha}:=(\alpha \otimes Id)\rho (v) \in A.
$$
It is well-known that:

(a) The coradical of $A$ is the linear span of the matrix
coefficients of
all simple $A$-comodules.

(b) The product in $A$ of two matrix coefficients is a matrix
coefficient of the tensor product. Specifically,
$$
\phi^V_{v, \alpha} \phi^W_{w, \beta} = \phi^{V\otimes
W}_{v\otimes w,
\alpha\otimes \beta}.
$$
It follows at once from (a) and (b) that $\Corad(A)$ is a
subalgebra of $A.$
Since the coradical is stable under the antipode, the claim
follows.

(3. $\Leftrightarrow$ 4.) To say that
$\Rad(H)$ is a Hopf ideal is equivalent to
saying that $\Corad(H^*)$ is a Hopf algebra, since
$\Corad(H^*)=(H/\Rad H)^*.$

(4. $\Rightarrow$ 1.) If $V,W$ are simple
$H-$modules then they factor
through
$H/\Rad(H).$ But $H/\Rad(H)$ is a Hopf algebra,
so $V\otimes W$ also factors through
$H/\Rad(H),$ so it is semisimple.

(3. $\Rightarrow$ 5.) Clear, since a cosemisimple Hopf
algebra is
involutory.

(5. $\Rightarrow$ 3.) Consider the subalgebra $B$ of $A$
generated by $\Corad(A).$
This is a Hopf algebra, and $S^2=Id$ on it.
Thus, $B$ is cosemisimple and hence
$B=\Corad(A)$ is a Hopf subalgebra of $A$.
\qed

\begin{Remark} {\rm The assumption that the base field has
characteristic $0$
is needed only in the proof of (5. $\Leftrightarrow$ 3.)
}
\end{Remark}

\section{The Classification of Triangular Hopf Algebras
With The Chevalley Property}

\subsection{The Main Theorem}

Our main result is the following theorem.

\begin{theorem}\label{main} Let $H$ be a
finite-dimensional triangular
Hopf algebra over $\C.$ Then the following are equivalent:
\ben
\item $H$ is a twist of a finite-dimensional triangular
Hopf algebra with $R-$matrix of
rank $\le 2$ (i.e. of a modified supergroup algebra).

\item $H$ has the Chevalley property.
\een
\end{theorem}

The proof of this theorem is contained in the next two
subsections.

\subsection{Local Finite-Dimensional Hopf Superalgebras
Are Exterior Algebras}

\begin{theorem}\label{1} Let $\cH$ be a local
finite-dimensional Hopf superalgebra (not necessarily
supercommutative). Then $\cH=\Lambda V^*$ for a finite-dimensional
vector space $V$. In other words, $\cH$ is the function algebra of
an odd vector space $V$.
\end{theorem}

\begin{remark} {\rm Note that in the commutative case
Theorem \ref{1}
is a special case of Proposition 3.2 of [Ko].
}
\end{remark}

\proof It is sufficient to show that $\cH^*=\Lambda V$ for some
vector space $V$, as $(\Lambda V)^*=\Lambda V^*$ as Hopf
superalgebras. For this, it is sufficient to show that $\cH^*$ is
generated by primitive elements, since the sub Hopf superalgebra
in $\cH^*$ generated by a basis of the space of primitive elements
of $\cH^*$ is clearly a {\it free} anti-commutative algebra on its
generators.

Let $I$ be the kernel of the counit in $\cH$. Then $I=\Rad(\cH)$
since $\cH$ is local. So in particular there exists a positive
integer $N$ such that for any $x_1,\dots,x_N\in \cH$ one has
$$(x_1-\eps(x_1)1)\cdots (x_N-\eps(x_N)1)=0.$$

Let $\delta_k: \cH^*\raro (\cH^*)^{\ot k}$ be the map dual to the
map $\cH^{\ot k}\raro \cH$ defined by $$x_1\ot \cdots \ot
x_k\mapsto (x_1-\eps(x_1)1)\cdots (x_k-\eps(x_k)1)$$ (this map was
introduced by Drinfeld in [Dr]). We see that we have a filtration
of $\cH^*:$ $\cH^*=\cup \cH^*_k,$ where $\cH^*_k$ is the kernel of
$\delta_k$ (the $N-$th term of this filtration is $\cH^*$). In
other words, $\cH_k^*$ is the orthogonal complement of $I^k.$

Let $V\subseteq \cH^*$ be the space of primitive elements, and
${\cal B}:=\Lambda V\subseteq \cH^*$ the corresponding Hopf
supersubalgebra generated by them. We will prove by induction in
$k$ that $\cH_k^*$ is contained in ${\cal B},$ which will complete
the proof.

The base of induction is obvious (as $\delta_1(x)=x-\eps(x),$
hence $\cH_1^*=\C$). Suppose the statement is known for $k=n,$ and
let $a\in \cH_{n+1}^*.$ Then it is straightforward to verify that
\linebreak $j:=\Delta(a)-a\otimes 1 - 1\otimes a\in \cH_n^*\otimes
\cH_n^*$. So by the induction assumption $j\in {\cal B}\otimes
{\cal B}$.

Since the second cohomology of the supercoalgebra $\Lambda V$ is
$S^2V$, any cohomology class can be represented by an even
cocycle. Thus, by subtracting from $a$ an element of $\Lambda V$,
if needed, we may assume that $a$ is even.

We claim that $a\in {\cal B}$. Suppose otherwise. Let $m\ge 2$ be
the smallest integer for which there exists a noncommutative
polynomial $Q(x)$ of elements of $\Lambda V$ and $x$, of degree
smaller than $m$ with respect to $x$, such that $a^m=Q(a)$ (i.e.
$Q(x)$ is a sum of expressions $b_1xb_2x\cdots b_{l-1}xb_l$,
$b_i\in {\cal B}$, $l\le m$); such polynomial exists because of
finite dimensionality.

Since $a$ is not in $\Lambda V$, there exists a linear functional
$f$ on the Hopf algebra which vanishes on $\Lambda V$ but equals
to $1/m$ on $a$. Then applying it in the second component of the
equation $\Delta(a^m)=\Delta(Q(a))$, we find that $a^{m-1}=P(a)$,
where $P$ is a noncommutative polynomial of degree $< m-1$. This
is a contradiction, hence $a\in {\cal B}$. We are done. \qed

\begin{remark}\label{ap}
{\rm In the appendix we give another proof of Theorem \ref{1}
using the Lifting method of [AS2].
}
\end{remark}

Theorem \ref{1} will be used in the next subsection, but
it also
allows one to give the following proof of Corollary
\ref{kostantf}.

\noin
{\bf Proof of Corollary \ref{kostantf}:}
Let $I$ be the ideal in $\cH^*$
generated by
all the odd elements. It is easy to see that this is a
Hopf ideal.
Consider the Hopf algebra $E:=\cH^*/I$ (the even part).
This is an {\em ordinary} commutative Hopf algebra, so
$E=F(G)$ for a
suitable
finite group $G.$ Moreover, it is clear that every
element of $I$
is nilpotent, so $I=\Rad(\cH^*).$ Thus, irreducible
$\cH^*-$modules are
$1-$dimensional, and are
parameterized by $g\in G.$ Let us call them $L_g.$ Also,
we see that $G={\bf G}(\cH).$

Let $P_g$ be the projective cover of the irreducible
module $L_g.$
Then $\cH^*=\bigoplus_g P_g,$ where $P_g$ are
indecomposable
two-sided ideals (the ideals are two-sided since the
algebra is commutative in the super-sense).
In particular, $P_g$ are local algebras, with $1-$dimensional
semisimple quotient. Also, we have a natural projection
of algebras
$\cH^*\raro P_g$ for all $g,$ in particular $\cH^*\raro
P_1.$

Note that $\cH$ acts on $\cH^*$ on the left and right. In
particular, so does the
group $G.$

\noin
{\bf Lemma.} {\em The following hold:
\ben
\item $g_1 P_g g_2=P_{g_1gg_2}.$
\item $\Delta(P_g)\subset
\bigoplus_{g_1,g_2:g_1g_2=g}P_{g_1}\otimes
P_{g_2}.$
\een
}

\noin
{\bf Proof:} Straightforward. \qed

\noin
{\bf Corollary.} {\em The ideal ${\mathcal
I}:=\bigoplus_{g\ne 1}P_g$
is a Hopf ideal, and
thus $P_1=\cH^*/{\mathcal I}$ is a Hopf superalgebra.
}

Thus, $P_1^*\subset \cH$ is a sub Hopf superalgebra
with an action of $G$, and we have a
factorization $\cH=\C[G]\ltimes P_1^*.$
The Hopf superalgebra $P_1$ is local, so $P_1^*=\Lambda V$ by
Theorem \ref{1}.
This concludes the proof. \qed

\begin{remark} {\rm
Here is the same proof,
described in a more intuitive geometric language.
Consider $\tilde
G:=\Spec(\cH^*).$ This is an
affine
supergroup scheme. Let $G\subseteq \tilde G$ be the even
part of
$\tilde G.$ Then $G$ is a finite group scheme, so by a
standard theorem it is
a finite group. Let $V$ be the
connected component of
the identity in $\tilde G.$ Then
the function algebra ${\cal O}(V)$ on $V$ is a local
finite-dimensional Hopf
superalgebra.
So by Theorem \ref{1}, ${\cal O}(V)=\Lambda V^*$ for some
finite-dimensional
vector space $V$.

Thus, we have a split exact sequence of algebraic supergroups

$$1\raro V\raro \tilde G\raro G\raro 1$$
(it is split because $G$ is a subgroup of $\tilde G$
complimentary to $V$).
So $\tilde G=G\ltimes V,$ as desired.
}
\end{remark}

\subsection{Proof of the Main Theorem}

We start by giving a super-analogue of Theorem 3.1 in [G].

\begin{lemma}\label{super3.1}
Let $\cH$ be a minimal triangular pointed Hopf
superalgebra.
Then $\Rad(\cH)$ is a Hopf ideal, and $\cH/\Rad(\cH)$ is
minimal triangular.
\end{lemma}

\proof The proof is a tautological generalization of the
proof of Theorem 3.1 in [G] to the super case.

First of all, it is clear that $\text{Rad}(\cH)$ is a Hopf
ideal,
since its orthogonal complement (the coradical of $\cH^*$)
is a sub Hopf superalgebra (as $\cH^*$ is isomorphic to
$\cH^{cop}$
as a coalgebra, and hence is pointed).
Thus, it remains to show that the triangular structure on
$\cH$
descends to a minimal triangular structure on
$\cH/\Rad(\cH)$.
For this, it suffices to prove that the
composition of the Hopf superalgebra maps
$$
\text{Corad}(\cH^{*cop})\hookrightarrow \cH^{*cop}\to \cH\to
\cH/\Rad(\cH)
$$
(where the middle map is given by the $R-$matrix) is an
isomorphism. But this follows from the fact that
for any surjective coalgebra map $\eta: C_1\to C_2$,
the image of the coradical of $C_1$ contains the coradical of
$C_2$ ([M], Corollary 5.3.5): One needs to apply this
statement to the map $\cH^{*cop}\to \cH/\text{Rad}(\cH)$.
\qed

\begin{lemma}\label{3} Let $\cH$ be a minimal
triangular pointed Hopf superalgebra, such that the \linebreak
$R-$matrix $\cR$ of $\cH$ is unipotent (i.e. $\cR-1\ot 1$ is $0$
in $\cH/\Rad(\cH)\ot \cH/\Rad(\cH)$). Then $\cH= \Lambda V$ as a
Hopf superalgebra, and $\cR=e^r,$ where $r\in S^2V$ is a
nondegenerate symmetric (in the usual sense) bilinear form on
$V^*$.
\end{lemma}

\proof By Lemma \ref{super3.1}, $\Rad(\cH)$ is a Hopf ideal, and
$\cH/\Rad(\cH)$ is minimal triangular. But the $R-$matrix of
$\cH/\Rad(\cH)$ must be $1\ot 1,$ so $\cH/\Rad(\cH)$ is
$1-$dimensional. Hence $\cH$ is local, so by Theorem \ref{1},
$\cH=\Lambda V.$ If $\cR$ is a triangular structure on $\cH$ then
it comes from an isomorphism $\Lambda V^*\raro \Lambda V$ of Hopf
superalgebras, which is induced by a linear isomorphism
$r:V^*\raro V.$ So $\cR=e^r,$ where $r$ is regarded as an element
of $V\ot V.$ Since $\cR \cR_{21}=1,$ we have $r+r^{21}=0$ (where
$r^{21}=-r^{op}$ is the opposite of $r$ in the supersense), so
$r\in S^2V$. \qed

\begin{remark}
{\rm The classification of pointed finite-dimensional Hopf
algebras with coradical of dimension $2$ is known [CD,N]. In the
appendix we use the Lifting method [AS1,AS2] to give an
alternative proof. Below we shall need the following more precise
version of this result in the triangular case. }
\end{remark}

\begin{lemma}\label{5} Let $H$ be a
minimal triangular pointed Hopf algebra,
whose coradical is $\C[\Z_2]=sp\{1,u\},$ where $u$ is the
Drinfeld element of $H.$
Then $H=\overline{(\Lambda V)^{{\cal J}}}$
with the triangular structure of Corollary \ref{mintr},
where
${\cal J}=e^{r/2},$ with $r\in S^2V$ a nondegenerate element.
In particular, $H$ is a twist of a modified supergroup algebra.
\end{lemma}

\proof Let $\cH$ be the associated triangular Hopf superalgebra
to $H,$ as described in Theorem \ref{cores}. Then the
$R-$matrix of $\cH$
is unipotent, because it turns into
$1\ot 1$ after killing the radical.

Let ${\cH}_m$ be the minimal part of $\cH.$ By
Lemma \ref{3}, ${\cH}_m=\Lambda V$ and $\cR=e^r,$ $r\in
S^2V.$ So if ${\cal J}:=e^{r/2}$ then
$\cH^{{\cal J}^{-1}}$ has $R-$matrix equal to
$1\ot 1.$ Thus, $\cH^{{\cal J}^{-1}}$ is cocommutative,
so by Corollary \ref{kostantf}, it equals
$\C[\Z_2]\ltimes \Lambda V.$ Hence
$\cH=\C[\Z_2]\ltimes (\Lambda V)^{{\cal J}}$,
and the result follows from Proposition \ref{1a}.
\qed

We shall need the following lemma.

\begin{lemma}\label{jh}
Let $B\subseteq A$ be finite-dimensional associative unital
algebras.
Then any simple $B-$module is a constituent
(in the Jordan-Holder series) of some simple $A-$module.
\end{lemma}

\proof Since $A,$ considered as a $B-$module, contains $B$
as a $B-$module, any simple $B-$module is a constituent of
$A.$

Decompose $A$ (in the Grothendieck group of $A$) into
simple $A-$modules: $A=\sum V_i.$ Further decomposing as
$B-$modules, we get $V_i=\sum W_{ij},$ and hence
$A=\sum_i\sum_j W_{ij}.$ Now, by Jordan-Holder theorem,
since $A$ (as a $B-$module) contains all simple
$B-$modules, any simple $B-$module $X$ is in $\{W_{ij}\}.$
Thus, $X$ is a constituent of some $V_i,$ as desired. \qed

\begin{proposition}\label{6} Any minimal
triangular Hopf algebra $H$ with the Chevalley property
is a twist of a triangular Hopf algebra with
$R-$matrix of rank $\le 2.$
\end{proposition}

\proof
By Proposition \ref{chev}, the coradical $H_0$ of $H$  is a Hopf
subalgebra,
since $H \simeq H^{*cop}$, being minimal triangular.
Consider the  Hopf algebra map $\varphi: H_{0} \to H^{*cop}/
\Rad(H^{*cop})$,
given by the composition of the following maps:
$$
H_{0} \hookrightarrow H \simeq H^{*cop} \to H^{*cop}/
\Rad(H^{*cop}),
$$
where the second map is given by the $R$-matrix. We claim that
$\phi$ is an
isomorphism. Indeed, $H_{0}$ and $H^{*cop}/ \Rad(H^{*cop})$ have the same
dimension, since
$\Rad(H^{*cop}) = (H_{0})^{\bot}$, and $\phi$ is injective, since
$H_{0}$
is semisimple by
[LR]. Let $\pi: H\to H_0$ be the associated projection.

We see, arguing exactly as in [G, Theorem 3.1], that  $H_0$ is
also minimal triangular, say with $R-$matrix $R_0$.

Now, by [EG1, Theorem 2.1], we can find a twist $J$ in $H_0\ot H_0$
such that $(H_0)^{J}$ is isomorphic to a group algebra
and has $R-$matrix $(R_0)^{J}$ of rank $\le
2$. Notice that here we are relying on Deligne's theorem, as mentioned in the
introduction.

Let us now consider $J$ as an element of $H_0\ot H_0$ and the twisted Hopf
algebra
$H^{J}$, which is again triangular.

The projection $\pi:H^{J} \to (H_0)^{J}$ is still a Hopf algebra map, and
sends
$R^J$ to $(R_0)^J.$ It induces a projection $(H^J)_m\to \C[\Z_2],$ whose
kernel
$K_m$ is contained in the kernel of $\pi$.
Because any simple $(H^J)_m-$module
is contained as a constituent in a simple $H-$module (see
Lemma \ref{jh}), $K_m=\Rad((H^J)_m).$ Hence,
$(H^J)_m$ is minimal triangular and $(H^J)_m/\Rad((H^J)_m)=(\C[\Z_2],
R_u).$ It follows, again by minimality, that $(H^J)_m$ is also pointed with
coradical
isomorphic to $\C[\Z_2]$.
So by Lemma \ref{5}, $(H^J)_m,$ and hence $H^J,$ can be
further twisted
into a triangular Hopf algebra with $R-$matrix of rank
$\le 2,$ as desired. \qed

Now we can prove the main theorem.

\noin
{\bf Proof of Theorem \ref{main}:}
(2. $\Rightarrow$ 1.) By Proposition
\ref{chev},
$H/\Rad(H)$ is a semisimple Hopf algebra. Let
$H_m$ be the minimal part
of $H,$ and $H_m'$ be the image of $H_m$ in
$H/\Rad(H).$ Then $H_m'$ is
a semisimple Hopf algebra.

Consider the kernel $K$ of the projection $H_m\to
H_m'.$
Then $K=\Rad(H)\cap H_m.$ This means that any
element $k\in K$ is
zero in any simple $H-$module. This implies that
$k$ acts by zero in any simple $H_m-$module,
since by Lemma \ref{jh}, any simple $H_m-$module
occurs as a constituent of some simple $H-$module. Thus,
$K$ is contained in $\Rad(H_m).$
On the other hand, $H_m/K$ is semisimple, so
$K=\Rad(H_m).$ This shows that $\Rad(H_m)$ is a Hopf
ideal. Thus, $H_m$ is
minimal triangular
satisfying the conditions of Proposition \ref{6}.
By Proposition \ref{6}, $H_m$ is a twist of a triangular Hopf algebra
with $R-$matrix of rank $\le 2$. Hence $H$ is a twist of a
triangular Hopf algebra with $R-$matrix of rank $\le 2$
(by the same twist), as desired.

(1. $\Rightarrow$ 2.) By assumption, $\Rep(H)$ is
equivalent to $\Rep(\tilde G)$ for some supergroup $\tilde G$
(as a tensor category without braiding).
But we know that supergroup algebras have the Chevalley
property, since, modulo their radicals, they are group
algebras.  This concludes the proof of the main
theorem. \qed

\begin{remark} {\rm Notice that it follows from the proof
of the
main
theorem
that any triangular Hopf algebra with the Chevalley
property can be
obtained by twisting of a triangular Hopf algebra
with $R-$matrix of rank $\le 2$ by an {\it even} twist.
}
\end{remark}

\begin{definition} If a triangular Hopf algebra
$H$ over $\C$ satisfies condition 1. or 2. of Theorem
\ref{main}, we will
say that $H$ is of supergroup type.
\end{definition}

\subsection{Corollaries of the Main Theorem}

\begin{corollary} \label{minim} A finite-dimensional triangular
Hopf
algebra $H$ is of supergroup type if and only if so is
its minimal part $H_m.$
\end{corollary}

\proof If $H$ is of supergroup type then
$\Rad(H)$ is a Hopf
ideal, so like
in the proof of Theorem \ref{main} (2. $\Rightarrow$ 1.) we
conclude
that $\Rad(H_m)$
is
a Hopf ideal, i.e. $H_m$ is of supergroup type.

Conversely, if $H_m$
is of supergroup type
then $H_m$ is a twist of a triangular Hopf algebra
with $R-$matrix of rank $\le 2$. Hence $H$ is a twist of a
triangular
Hopf algebra with $R-$matrix of rank $\le 2$
(by the same twist), so $H$ is
of supergroup type. \qed

\begin{corollary}\label{6a} A finite-dimensional
triangular Hopf algebra
whose coradical is a Hopf subalgebra is of supergroup type.
In particular, this is the case for finite-dimensional triangular
pointed Hopf algebras.
\end{corollary}

\proof This follows from Corollary \ref{minim}.\qed

\begin{corollary} \label{triangcopoint} Any
finite-dimensional triangular
basic Hopf algebra is of supergroup type.
\end{corollary}

\proof A basic Hopf algebra automatically has
the Chevalley property
since all its irreducible modules are
$1-$dimensional. Hence the result follows from the
main theorem. \qed

\subsection{Questions}

The above results motivate the following question.

\begin{question} \label{q1}
Does any finite-dimensional
triangular Hopf algebra over $\C$ have the Chevalley property
(i.e. is of supergroup type)?
Is it true under the assumption that $S^4=Id$
or at least under the assumption that $u^2=1$?
\end{question}

\begin{remark} {\rm Recall [G] that it is not known
whether any finite-dimensional triangular Hopf algebra over $\C$
has the property $u^2=1$ or at least  $S^4=Id$.
It is also not known if $S^4=Id$ implies $u^2=1$ for
triangular Hopf algebras. However, it is clear that for
finite-dimensional triangular Hopf algebras $H$ of supergroup type,
$u^2=1$ (and hence $S^4=Id$). Indeed, since $S^2=Id$ on the semisimple
part of $H,$ $u$ acts by a scalar in any irreducible representation of
$H.$ In fact, since $\tr(u)=\tr(u^{-1}),$ we have that $u=1$ or $u=-1$ on
any irreducible representation of $H,$ and hence $u^2=1$ on any irreducible
representation of $H.$ Thus, $u^2$ is unipotent. But it is of finite order
(as it is a grouplike element), so it is equal to $1$ as desired.
}
\end{remark}

\begin{remark} {\rm Note that the answer to question \ref{q1}
is negative
in the infinite-dimensional case. Namely, although the
answer is
positive
in the cocommutative case (by [C]), it is negative already
for triangular Hopf algebras with $R-$matrix of rank $2,$
which
correspond to cocommutative Hopf superalgebras.
Indeed, let us take the cocommutative Hopf superalgebra
$\cH:=\text{U}(\text{gl}(n|n))$
(for the definition of the Lie superalgebra
$\text{gl}(n|n),$ see [Ka, p.29]). The associated
triangular Hopf algebra $\ol {\cH}$ does not have the
Chevalley property,
since it is well known that Chevalley theorem fails
for Lie superalgebras (e.g. $\text{gl}(n|n)$); more
precisely, already
the product of the vector and covector representations
for this Lie superalgebra is not semisimple.
}
\end{remark}

\begin{remark} {\rm It follows from Corollary \ref{minim}
that a positive answer to Question \ref{q1} in the minimal case
would imply the general positive answer.
}
\end{remark}

Here is a generalization of Question \ref{q1}.

\begin{question}\label{q2} {\rm Does any $\C-$linear abelian
symmetric rigid tensor category, with $\text{End}({\bf 1})=\C$
and finitely many
simple objects, have the Chevalley property? }
\end{question}

Even a more ambitious question:

\begin{question}\label{q2a}{\rm
Is such a category equivalent to the category of
representations of a
finite-dimensional triangular
Hopf algebra with $R-$matrix of rank $\le 2$? In particular,
is it equivalent to the category of representations of a
supergroup, as a category without braiding?
Are these statements valid at least for categories
with Chevalley property?
For semisimple categories?}
\end{question}

\section{Finite-Dimensional Cotriangular Pointed Hopf Algebras
are Generated by Grouplikes and Skewprimitives}

There is a conjecture (see [AS2]) that any
finite-dimensional pointed Hopf algebra over $\C$ is
generated by
grouplike and skew primitive elements. Here we confirm it
in the
cotriangular case.

\begin{Theorem}\label{gen} A finite-dimensional
cotriangular pointed Hopf algebra
$H$ over $\C$ is generated by grouplike and skew primitive
elements.
\end{Theorem}

In order to prove the theorem, we will need the following
lemma.

\begin{Lemma}\label{l1}
Let $H$ be a finite-dimensional pointed Hopf
algebra or superalgebra. Then the
following are equivalent:
\ben
\item $H$ is generated by grouplike and skew
primitive elements.
\item There exists a faithful $H^*-$module which
is a direct sum of
tensor products of $H^*-$modules of dimension $2.$
\een
\end{Lemma}
\proof Irreducible $H^*-$modules are
$1-$dimensional, so a
$2-$dimensional representation has the form
$$ a\mapsto
\left(
\begin{array}{ll}
p(a) & r(a)\\
0 & q(a)
\end{array}
\right ),\;
a\in H^*,
$$
where $p,q$ are characters (i.e. belong to ${\bf
G}(H)$), and
$r$ is a $q:p$ skew
primitive. Conversely, such a $2-$dimensional representation
exists for any skew primitive element.
Matrix elements
of tensor products of representations of $H^*$
are products of the
matrix elements
of these representations (as elements of $H$).
This implies the lemma. \qed

Now we are ready to give a proof of Theorem \ref{gen}.

We need to show that $H^*$ has a faithful
representation which
is a direct sum of products of $2-$dimensional
representations. By [G], the Drinfeld element
$u$ of $H^*$ satisfies $u^2=1.$ Hence
by Theorem \ref{cores0}, we can
replace $H^*$ by the corresponding Hopf superalgebra
$\widetilde {H^{*}}$ (as this does not change the
representation category as a tensor category).
But $H^*$ is basic triangular, which means,
by Corollary \ref{triangcopoint}, that
$\widetilde {H^{*}}$ is twist equivalent to a supergroup
algebra $B.$
Thus by Lemma \ref{l1}, it suffices to show that $B^*$
(the dual of $B$) is
generated by
grouplikes and skew primitives.

But $B=\C[G]\ltimes
\Lambda V$ where
$G$ is abelian.
Thus, $V$ is decomposed in the direct sum of
eigenspaces for $G.$
So let $v_1,...,v_n$ be a basis of $V$, such that
$gv_ig^{-1}=\chi_i(g)v_i$, where $\chi_i$ are some characters
of $G$. Using this presentation of $B$, it is easy to compute
its dual $B^*$ and show that it is generated as an
algebra by $G^\vee$ (the character
group) and $\chi_i:1$ skew primitive elements $\xi_i$,
$i=1,...,n$.
We are done. \qed

\begin{Corollary}\label{mtp}
Theorem 5.1 of [G] gives
the classification of all minimal triangular pointed
Hopf algebras.
\end{Corollary}

\proof
Since minimal triangular pointed Hopf algebras
are also cotriangular, by Theorem \ref{gen},
they are generated by grouplikes and skewprimitives
(which answers a question from [G]).
On the other hand, [G, Theorem 5.1]
gives a classification of minimal triangular Hopf algebras
which are generated by grouplikes and skewprimitives.
\qed

\begin{Remark} {\rm Lemma \ref{l1} implies
that if $H_1,H_2$ are finite-dimensional pointed
Hopf algebras, and
the comultiplication of $H_1^*$ is obtained by
conjugating that of $H_2^*$
by an invertible element (not necessarily a
twist), then
$H_1$ is generated by
grouplike and skew primitive elements if and only
if so is $H_2.$
}
\end{Remark}

\section{Categorical Dimensions in Symmetric
Categories\\ With Finitely Many Irreducibles are Integers}
In this paper we classified finite-dimensional
triangular Hopf algebras
with the Chevalley property.
In conclusion, let us give one result that is valid in
a greater generality: for any finite-dimensional
triangular Hopf
algebra, and even for any symmetric rigid category
with finitely many irreducible objects.

Let ${\cal C}$ be a $\C-$linear abelian symmetric rigid
category with ${\bold 1}$ as its unit
object, and suppose that $\End({\bold 1})=\C.$ Recall
that there is a natural
notion of dimension in ${\cal C}$, generalizing the
ordinary dimension of an object
in $\Vect.$ Let $\beta$ denote the
commutativity constraint in ${\cal
C},$ and for an object
$V,$ let $ev_V,$ $coev_V$ denote the associated
evaluation and coevaluation
morphisms.

\begin{Definition}\label{dimcat} {\bf [DM]}
The categorical dimension $\text
{dim}_c(V)\in \C$ of $V\in \text{Ob}({\cal C})$ is the
morphism
\begin{equation}\label{dimc}
\text {dim}_c(V):{\bold 1}\stackrel{ev_V}{\longrightarrow}V\ot
V^*\stackrel{\beta_{V,V^*}}{\longrightarrow}V^*\ot
V\stackrel{coev_V}{\longrightarrow}{{\bold 1}}.
\end{equation}
\end{Definition}

The main result of this section is the following:

\begin{Theorem}\label{cd} In any $\C-$linear abelian symmetric
rigid tensor category ${\cal C}$
with finitely
many irreducible objects, the categorical dimensions of
objects are
integers.
\end{Theorem}

\proof First note that the categorical dimension of any
object $V$ of ${\cal C}$ is an
algebraic integer. Indeed, let $V_1\dots,V_n$ be
the irreducible objects of ${\cal
C}.$ Then $\{V_1\dots,V_n\}$ is a basis of the
Grothendieck ring of ${\cal C}.$
Write $V\otimes V_i=\sum_j N_{ij}(V)V_j$ in the
Grothendieck ring. Then $N_{ij}(V)$ is
a matrix with integer entries, and $\dim_c(V)$ is
an eigenvalue of this matrix.
Thus, $\dim_c(V)$ is an algebraic integer.

Now, if $\dim_c(V)=d$ then it is easy to show (see e.g.
[De2]) that
$$\dim_c(S^kV)=d(d+1)\cdots (d+k-1)/k!,$$
and
$$\dim_c(\Lambda^k V)=d(d-1)\cdots (d-k+1)/k!,$$
hence they are also algebraic
integers. So the theorem follows from:

\noin
{\bf Lemma.} {\em Suppose $d$ is an algebraic
integer such that $d(d+1)\cdots (d+k-1)/k!$
and \linebreak $d(d-1)\cdots (d-k+1)/k!$ are algebraic
integers for all $k.$ Then $d$ is an
integer.}

\noin
{\bf Proof:} Let $Q$ be the minimal monic polynomial of
$d$ over $\Z.$
Since $d(d-1)\cdots (d-k+1)/k!$ is an algebraic
integer, so are the
numbers $d'(d'-1)\cdots (d'-k+1)/k!,$ where $d'$ is
any algebraic conjugate of $d$.
Taking the product over all conjugates, we get that
$$N(d)N(d-1)\cdots N(d-k+1)/(k!)^n$$
is an integer,
where $n$ is the degree of
$Q.$
But $N(d-x)=(-1)^nQ(x).$ So we get that
$Q(0)Q(1)\cdots Q(k-1)/(k!)^n$ is an integer.
Similarly from the identity for $S^kV,$ it
follows that $Q(0)Q(-1)\cdots Q(1-k)/(k!)^n$ is an
integer. Now, without loss of generality, we can
assume that
$Q(x)=x^n+ax^{n-1}+...,$ where $a\le 0$
(otherwise replace $Q(x)$ by
$Q(-x);$ we can do it since our condition is
symmetric under this
change).
Then for large $k,$ we have $Q(k-1)<k^n,$ so the
sequence
$b_k:=Q(0)Q(1)\cdots Q(k-1)/k!^n$ is decreasing in
absolute value
or zero starting from some place. But a sequence
of integers cannot
be strictly decreasing in absolute value forever.
So $b_k=0$ for some
$k,$ hence $Q$ has an integer root. This means
that $d$ is an integer
(i.e. $Q$ is linear), since $Q$ must be irreducible
over the rationals.
This concludes the proof of the lemma, and hence
of the theorem. \qed

\begin{Corollary}\label{cdtr} For any triangular Hopf
algebra $H$ (not
necessarily
finite-dimensional), the categorical dimensions of
its finite-dimensional
representations are integers.
\end{Corollary}

\proof It is enough to consider the minimal part $H_m$ of
$H$ which is
finite-dimensional, since $\dim_c(V)=\tr(u_{|V})$ for any
module $V$ (where $u$ is the Drinfeld element of $H$), and
$u\in H_m.$
Hence the result follows from Theorem \ref{cd}.
\qed

\begin{Remark} {\rm Theorem \ref{cd} is false without the
finiteness conditions. In fact, in this case any complex
number can be a dimension, as is
demonstrated in examples constructed by Deligne [De1,
p.324-325]. Also, it is well known that the theorem
is false for ribbon, nonsymmetric categories
(e.g. for fusion categories of
semisimple representations of finite-dimensional quantum
groups at roots of unity,
where dimensions can be irrational algebraic integers).
}
\end{Remark}

\begin{Remark} {\rm Note that Theorem \ref{cd}
can be regarded as a piece of supporting evidence for a
positive
answer to
Question \ref{q2a}.
}
\end{Remark}

\begin{Remark} {\rm In any rigid braided tensor category with
finitely many
irreducible objects, one can define the Frobenius-Perron
dimension
of an object $V$, $\FPdim(V),$ to be the largest positive
eigenvalue of the
matrix of multiplication by $V$ in the Grothendieck ring.
This dimension is well defined
by the Frobenius-Perron theorem, and has the usual
additivity and
multiplicativity properties.
For example, for the category of representations of a
quasi-Hopf algebra, it is just the usual dimension of the
underlying vector space.
If the answer to
Question \ref{q2a} is positive then $\FPdim(V)$
for symmetric categories is always an integer,
which is equal to $\dim_c(V)$ modulo $2.$ It would be
interesting
to check this, at least in the case of modules over a
triangular
Hopf algebras, when the integrality of $\FPdim$ is automatic
(so only the mod $2$ congruence has to be checked).}
\end{Remark}

\section{Appendix: On Pointed Hopf Algebras}

In this appendix we use the Lifting method [AS1,AS2] to
give other proofs of Theorem \ref {1} and Corollary \ref
{mtp}, and a generalization of Lemma \ref {5}.

By {\em a braided Hopf algebra} we shall mean a Hopf algebra in
the braided tensor category of Yetter-Drinfeld modules over a
group algebra $\C[\Gamma]$, where $\Gamma$ is a finite abelian
group. For example, one can endow the exterior algebra $\Lambda V$
with the structure of a braided Hopf algebra, as follows. Let
$x_1, \dots, x_N$ be a basis of $V$ and let there be given $g_1,
\dots, g_N\in \Gamma$ and $\chi_1, \dots, \chi_N\in \Gamma^\vee$
such that
$$
\chi_i(g_j) = -1,\quad 1\le i, j\le N.
$$
Then $V$ is a Yetter-Drinfeld module over $\C[\Gamma]$
where the action and coaction of $\Gamma$ on $x_i$
are given by $\chi_i$ and $g_i$ respectively.
These action and coaction extend to $\Lambda V$, and
turn $\Lambda V$ into a braided Hopf algebra.

\begin{Lemma}\label{as1} Let $R= \bigoplus_{n\ge 0} R(n)$
be a graded braided Hopf algebra
such that $R(0) = \C$, $R(1) \simeq V$ as a
Yetter-Drinfeld module (with
the assumptions above), and $R$ is generated by $R(1)$.
Then $R$ is
isomorphic to  $\Lambda V$ as a graded braided Hopf
algebra.
\end{Lemma}

\proof It is known, and easy to see, that $\Lambda V$
satisfies all the hypotheses that
$R$ does, plus that the primitive elements are
concentrated in degree one:
${\bf P}(\Lambda V) = \Lambda V(1) = V$ (see for instance
[AS1, Section 3]). In other words,
$\Lambda V$ is the Nichols algebra of $V,$ and there
exists a surjective homomorphism
of graded braided Hopf algebras $R\to \Lambda V$ which
is the identity in
degree one (see for instance [AS2, Lemma 5.5]). On
the other hand, it is
clear that $\Lambda V$ can be presented by generators
$x_1, \dots, x_N$ with relations
\begin{equation}\label{relslambdaV}
x_ix_j + x_jx_i =0,\;1\le i, j\le N.
\end{equation}
So in particular $x_i^2 = 0$ for all $i$. It is thus
enough to show that
equations (\ref{relslambdaV}) also hold in $R$, with an
evident abuse of
notation. But $x_ix_j + x_jx_i$
is a primitive element of $R$, whose action is given by
the character
$\chi_i\chi_j,$ and whose coaction is given by $g_ig_j$.
Since
$\chi_i\chi_j(g_ig_j) = 1$ and $R$ is finite-dimensional,
$x_ix_j + x_jx_i = 0$ in $R$ by [AS1, Lemma 3.1]. \qed

Let $H$ be a finite-dimensional pointed Hopf algebra such
that ${\bf G}(H)$ is isomorphic to $\Gamma$.
We recall that the Lifting method [AS1, AS2] attaches to
$H$ several invariants:

\begin{itemize}
\item The graded Hopf algebra $\gr H$, associated to the
coradical filtration of $H$.

\item A graded braided Hopf algebra $R;$ the coinvariants
of the homogeneous projection from $\gr H$ to
$\C[\Gamma]$.

\item A Yetter-Drinfeld module $W:=R(1)$ over $\C
[\Gamma]$, called the
infinitesimal braided vector space of $H$.
\end{itemize}

We conclude immediately from Lemma \ref{as1}:

\begin{Corollary}\label{as2} Let $H$ be a finite-dimensional
pointed Hopf
algebra such that ${\bf G}(H)$ is isomorphic to $\Gamma$.
Assume that the
infinitesimal braiding of $H$ is isomorphic to $V$ as
above. Then $H$ is
generated by grouplike and skewprimitive elements. \qed
\end{Corollary}

\begin{Remark} {\rm Notice that Corollary \ref{as2}
allows one to give an alternative proof of Corollary \ref{mtp}.
For, Lemmas 5.3 and 5.4 in [G] imply that the infinitesimal
braiding of any minimal triangular pointed Hopf algebra is
isomorphic to a $V$ as above. }
\end{Remark}

Assume now that $\Gamma =\Z_2$. There is only one
possible choice for $V$
as above, namely
$g_1 = \dots = g_N = u$ and $\chi_1 = \dots =\chi_N = $
the sign. This gives
the Hopf superalgebra as explained in Section 5.
Let now $H$ be a
finite-dimensional pointed Hopf algebra such that
${\bf G}(H)$
is isomorphic to
$\Z_2$. Then the infinitesimal braiding of $H$ is
isomorphic to  $V$ as
above by [AS1, Lemma 3.1] again, for some natural
number $N$.
The Lifting method gives a very direct proof of the
following well-known result.

\begin{Theorem}\label{as3} {\bf [N, Th. 4.2.1], [CD]}
If
$H$ is a finite-dimensional pointed Hopf algebra with
${\bf G}(H)$ isomorphic to
$\Z_2$, then $H \simeq \C[\Z_2]\ltimes \Lambda V$.
\end{Theorem}

\proof By the above remarks and Corollary \ref{as2},
we know that $\gr H
\simeq \C[\Z_2]\ltimes \Lambda V$ for some $V$. The fact
that
$H \simeq \gr H$
("there are no liftings" in the jargon of the Lifting
method) is a
particular case of the main result [AS1, Th.
5.5]. \qed

We can now give another proof of Theorem \ref{1}.

 It is enough to show that ${\cal H}^*=\Lambda V$ for
some $V$ as above, since $(\Lambda V)^*=\Lambda V^*$ as Hopf
superalgebras. By the hypothesis, the coradical of ${\cal H}^*$ is
trivial: $\Corad({\cal H}^*) = \C1$. We can consider the
biproduct $H:= \C[\Z_2]\ltimes {\cal H}^*$; that is, $H=\ol
{\cH}$ in our notation. We claim that $H$ is a finite-dimensional
pointed Hopf algebra with ${\bf G}(H)$ isomorphic to $\Z_2$.
Indeed, the filtration
$$\C [\Z_2]\subset
\C[\Z_2]\ltimes ({\cal H}^*)_1\subset\cdots
\C[\Z_2]\ltimes ({\cal H}^*)_j \cdots$$
is a coalgebra filtration, where $1 \subset ({\cal H}^*)_1
\subset\cdots  ({\cal H}^*)_j
\cdots$ is the coradical filtration of ${\cal H}^*$. Hence $\C
[\Z_2]$ contains
the coradical of $H$ and the other inclusion is evident.

It follows then from Theorem \ref{as3} that $H\simeq
\C[\Z_2]\ltimes \Lambda V$.
By [AS2, Lemma 6.2], ${\cal H}^* \simeq\Lambda V$ as braided
Hopf algebras, that is as Hopf superalgebras. \qed

\end{document}